\documentclass[10pt]{amsart}
\usepackage{amscd}
\usepackage{amssymb}
\usepackage{latexsym}
\usepackage{graphicx}
\usepackage[parfill]{parskip} 

\newcommand{\FT}{\,\, \widehat{} \,\,}
\newcommand{\init}{\vert_{t = 0}}

\newcommand{\abs}[1]{\left\vert #1 \right\vert}
\newcommand{\fixedabs}[1]{\vert #1 \vert}

\newcommand{\norm}[1]{\left\Vert #1 \right\Vert}
\newcommand{\fixednorm}[1]{\Vert #1 \Vert}
\newcommand{\bignorm}[1]{\bigl\Vert #1 \bigr\Vert}

\newcommand{\C}{\mathbb{C}}
\newcommand{\R}{\mathbb{R}}

\newcommand{\innerprod}[2]{\left\langle \, #1 , #2 \,
\right\rangle}

\newcommand{\biginnerprod}[2]{\bigl\langle \, #1 , #2
\, \bigr\rangle}

\newcommand{\angles}[1]{\langle #1 \rangle}

\DeclareMathOperator{\re}{Re}

\newtheorem{theorem}{Theorem}
\newtheorem{lemma}{Lemma}
\newtheorem{corollary}{Corollary}

\theoremstyle{definition}

\theoremstyle{remark}
\newtheorem{remark}{Remark}

\title[DKG 1d global well-posedness]{Global well-posedness below the charge norm for the Dirac-Klein-Gordon system in one space dimension}

\author{Sigmund Selberg}
\address{Department of Mathematical Sciences\\ Norwegian University of Science and Technology\\ Alfred Getz' vei 1\\ N-7491 Trondheim\\ Norway}
\email{sselberg@math.ntnu.no}
\urladdr{www.math.ntnu.no/~sselberg}
\thanks{Supported by Research Council of Norway, project 160192/V30, PDE and Harmonic Analysis}


\begin{document}

\begin{abstract}
We prove global well-posedness below the charge norm (i.e., the $L^2$ norm of the Dirac spinor) for the Dirac-Klein-Gordon system of equations (DKG) in one space dimension. Adapting a method due to Bourgain, we split off the high frequency part of the initial data for the spinor, and exploit nonlinear smoothing effects to control the evolution of the high frequency part. To prove the nonlinear smoothing we rely on the null structure of the DKG system, and bilinear estimates in Bourgain-Klainerman-Machedon spaces.
\end{abstract}

\maketitle

\section{Introduction}\label{Section1}

We study the Dirac-Klein-Gordon system (DKG) in one space dimension,
\begin{equation}\label{DKG1}
\left\{
\begin{alignedat}{2}
&D_t \psi + \alpha D_x  \psi + M\psi= \phi \beta\psi,&
\qquad\qquad
&\left(D_t=-i \partial_t, \,\, D_x=-i \partial_x\right)
\\
&- \square \phi + m^2 \phi= \innerprod{\beta\psi}{\psi}_{\C^2},
&
&\left(\square = -\partial_t^2 + \partial_x^2 \right)
\end{alignedat}
\right.
\end{equation}
with initial data
\begin{equation}\label{data}
  \psi \init = \psi_0 \in H^{-s}, \qquad \phi \init =
\phi_0 \in H^r, \qquad
\partial_t \phi \init = \phi_1 \in H^{r-1},
\end{equation}
where $r,s \ge 0$ are fixed later.
Here $\phi(t,x)$ is real-valued and $\psi(t,x) \in \C^2$ is the Dirac spinor, regarded as a column vector with components $\psi_1$, $\psi_2$. Further, $M,m \ge 0$ are constants, $H^r = (I-\partial_x^2)^{-r/2}L^2$ is the standard Sobolev space of order $r$ on $\R$, $\innerprod{\cdot}{\cdot}_{\C^2}$ is the standard inner product on $\C^2$, and we use the following representation of the Dirac matrices:
$$
  \alpha =    \begin{pmatrix}
    0 & 1  \\
    1 & 0
  \end{pmatrix},
  \qquad
  \beta =  \begin{pmatrix}
    1 & 0  \\
    0 & -1
  \end{pmatrix}.
$$

Global well-posedness (GWP) for DKG in 1d was first proved by Chadam \cite{Chadam:1973}, for data
$$
  (\psi_0,\phi_0,\phi_1) \in H^1 \times H^1 \times L^2.
$$
In recent years, many authors have improved this result, in the sense that the required regularity on the data has been lowered: see Table \ref{Table0} for an overview (here and elsewhere we use the convention that if the regularity of $\phi_1$ is not specified, it is understood to be always one order less than that of $\phi_0$).

\begin{table}
\caption{Global well-posedness for DKG in 1d}
\label{Table0}
\def\arraystretch{1.3}
\begin{center}
\begin{tabular}{|r|c|l|}
\hline
& $\psi_0 \in$ & $\phi_0 \in$
\\
\hline\hline
Chadam \cite{Chadam:1973}, 1973 & $H^1$ & $H^1$
\\
\hline
Bournaveas \cite{Bournaveas:2000}, 2000 & $L^2$ & $H^1$
\\
\hline
Fang \cite{Fang:2004b}, 2004 & $L^2$ & $H^r$, $1/2 < r \le 1$
\\
\hline
Bournaveas and Gibbeson \cite{BournaveasGibbeson:2006}, 2006 & $L^2$ & $H^r$, $1/4 \le r \le 1$
\\
\hline
Machihara \cite{Machihara:2006}, Pecher \cite{Pecher:2006}, 2006 & $L^2$ & $H^r$, $0 < r \le 1$
\\
\hline
\end{tabular}
\end{center}
\end{table}

All these GWP results rely on the conservation of charge:
\begin{equation}\label{conservation_of_charge}
  \norm{\psi(t)}_{L^2} = \norm{\psi(0)}_{L^2},
\end{equation}
which holds for smooth solutions decaying sufficiently fast at spatial infinity, hence also for any solution evolving from initial data for which local well-posedness (LWP) holds, provided the regularity of the spinor is at least $L^2$, of course.

By combining \eqref{conservation_of_charge} with a sufficiently strong a priori estimate for $\phi$, one can reduce the problem of GWP to proving LWP, again provided the regularity of the spinor is at least $L^2$. In fact, by the energy inequality for the wave equation,
\begin{equation}\label{phi_energy_estimate}
  \norm{\phi(t)}_{H^r}
  + \norm{\partial_t \phi(t)}_{H^{r-1}}
  \le
  C_T \left( \norm{\phi_0}_{H^r}
  + \norm{\phi_1}_{H^{r-1}}
  + \int_0^t \norm{\square \phi(\sigma)}_{H^{r-1}} \, d\sigma \right),
\end{equation}
for $0 \le t \le T < \infty$, assuming existence up to time $T > 0$. If $r < 1/2$, one can then apply the product law for Sobolev spaces (in one space dimension):
\begin{equation}\label{Sobolev_product}
  \norm{fg}_{H^{-c}} \le C_{a,b,c} \norm{f}_{H^a} \norm{g}_{H^b}
  \qquad \text{if}
  \qquad
  \left\{
  \begin{aligned}
  &a+b+c > 1/2,
  \\
  &a+b \ge 0,
  \quad b+c \ge 0,
  \quad a+c \ge 0.
  \end{aligned}
  \right.
\end{equation}
Thus, if $r < 1/2$, combining  \eqref{Sobolev_product}, \eqref{phi_energy_estimate} and \eqref{conservation_of_charge} gives
\begin{equation}\label{phi_apriori}
\norm{\phi(t)}_{H^r}
  + \norm{\partial_t \phi(t)}_{H^{r-1}}
  \le
  C_T \left( \norm{\phi_0}_{H^r}
  + \norm{\phi_1}_{H^{r-1}}
  + T \norm{\psi_0}_{L^2}^2 \ \right),
\end{equation}
for $0 \le t \le T$. This together with the conservation of charge \eqref{conservation_of_charge} shows that if LWP holds with data $\psi_0 \in L^2$ and $\phi_0 \in H^r$, then GWP follows immediately, if $r < 1/2$. The same is in fact true if $1/2 \le r \le 1$; we remark on this below.

Concerning LWP, the best result is the following. Note that we use $-s$ to denote the order of the Sobolev space for the spinor, since we are primarily interested in spaces of negative order.

\begin{theorem}\label{LWP_theorem}
The DKG system \eqref{DKG1} is LWP for data
$\,\psi_0 \in H^{-s}$, $\phi_0 \in H^r$ provided
$$
  s < \frac{1}{4}, \qquad r>0,
  \qquad
  \abs{s}
  \le r \le 1-s.
$$
\end{theorem}

In the most interesting case $s \ge 0$ the attributions are as follows:
Bournaveas \cite{Bournaveas:2000} proved LWP for $s = 0$ and $r=1$;
Fang \cite{Fang:2004b} proved the result for $0 \le s < 1/4$ and $1/2 < r \le 1-2s$;
Bournaveas and Gibbeson \cite{BournaveasGibbeson:2006} obtained $s = 0$ and $1/4 \le r \le 1$;
Machihara \cite{Machihara:2006} proved LWP for $0 \le s < 1/4$ and  $2\abs{s} \le r \le 1 - 2s$;
Pecher \cite{Pecher:2006}, independently of Machihara, obtained $0 \le s < 1/4$ and $\abs{s} \le r < 1 - 2s$; finally, it was shown by A.\ Tesfahun and the present author in \cite{Selberg:2006e} that the upper bound on $r$ in Pecher's result can be relaxed to $r \le 1-s$. Moreover, it was shown in \cite{Selberg:2006e} that the conditions in Theorem \ref{LWP_theorem} are optimal (up to, possibly, the endpoint $r=s=0$) if one iterates in Bourgain-Klainerman-Machedon spaces.

As remarked already, LWP implies GWP when the regularity of the spinor is at least $L^2$, i.e., when $s \le 0$ in Theorem \ref{LWP_theorem}, essentially due to the conservation of charge.

\begin{theorem}\label{GWP_theorem} (Cf.\ Table \ref{Table0}.)
The DKG system \eqref{DKG1} is GWP for data
$\,\psi_0 \in H^{-s}$, $\phi_0 \in H^r$ if
$$
  s \le 0, \qquad r>0,
  \qquad
  \abs{s}
  \le r \le 1-s.
$$
\end{theorem}

\begin{remark}
If $s=0$, this follows directly from Theorem \ref{LWP_theorem} if $0 < r < 1/2$, by the argument shown above, and the case $1/2 \le r \le 1$ then reduces to the case $0 < r < 1/2$ by using the propagation of higher regularity, which is discussed in the next remark. The case $s < 0$ can be reduced to $s=0$, again by propagation of higher regularity.
\end{remark}

\begin{remark}\label{Remark2} Note that the conditions in Theorem \ref{LWP_theorem} describe a convex region $R$ of the $(s,r)$-plane. Suppose that $(s,r)$ and $(s',r')$ both belong to $R$, and that $-s' \le -s$ and $r' \le r$. Assume that $\psi_0 \in H^{-s}$ and $\phi_0 \in H^r$. Then $\psi_0 \in H^{-s'}$ and $\phi_0 \in H^{r'}$ also, of course. Now suppose that the solution in the latter data space exists up to some time $0 < T' < \infty$, so
$$
  \psi \in C\bigl([0,T'];H^{-s'}\bigr), \qquad
  (\phi,\partial_t \phi) \in  C\bigl([0,T'];H^{r'} \times H^{r'-1}\bigr).
$$
Then in fact the solution retains the higher regularity of the initial data throughout the time interval $[0,T']$, so that
$$
  \psi \in C\bigl([0,T'];H^{-s}\bigr), \qquad
  (\phi,\partial_t \phi) \in  C\bigl([0,T'];H^{r} \times H^{r-1}\bigr).
$$
This property, known as propagation of higher regularity, is included in the LWP result of Theorem \ref{LWP_theorem}, since it is proved by an iteration scheme using estimates on the iterates to get a contraction.
\end{remark}

In this paper, we shall extend Theorem \ref{GWP_theorem} to a range of negative order Sobolev exponents for the spinor. 

Our main result is the following.

\begin{theorem}\label{Main_theorem}
The DKG system \eqref{DKG1} is GWP for data
$\,\psi_0 \in H^{-s}$, $\phi_0 \in H^r$ if
$$
  0 < s < \frac{1}{8},
  \qquad
  s + \sqrt{s^2+s} < r \le 1-s.
$$
\end{theorem}

Since $\psi_0$ is in a Sobolev space of negative order, the conservation law \eqref{conservation_of_charge} is not directly applicable. Bourgain \cite{Bourgain:1998} devised a method for dealing with situations where a PDE has a conservation law, but the initial data are too rough to use it directly. The essence of Bourgain's idea is to split the PDE into two parts, a `good' part and a `bad' part, by separating the low and high frequency parts of the data. For the good part one can use the conservation law, while for the `bad' part there must be a nonlinear smoothing effect in order for the method to work, so that the smoothed out nonlinear part can be fed back into the `good' part at the end of the possibly short time interval of existence of the `bad' part. Then this procedure is iterated, to reach a given, final time. This is the underlying idea, but in order to make it work, i.e., to close the finite induction, one needs sufficiently good estimates, of course.

The implementation of Bourgain's method for DKG requires an additional idea, since there is no conservation law for the field $\phi$, only for the spinor $\psi$. We overcome this difficulty by using the estimate \eqref{phi_apriori} for the inhomogeneous part, together with an additional induction argument involving a cascade of free waves; see Section \ref{Conclusion}.

To prove the nonlinear smoothing we rely on the null structure of the DKG system, which was recently completed by P.\ D'Ancona, D.\ Foschi and the present author in \cite{Selberg:2006b}, and on bilinear estimates in Bourgain-Klainerman-Machedon spaces.

The rest of this paper is organized as follows: In Section \ref{Preliminaries} we split the system into two separate systems, by a frequency cut-off on the spinor data. Then in Sections \ref{uw_estimates} and \ref{vz_estimates} we prove the main estimates for the split DKG system. In Section \ref{Conclusion} we put everything together to finish the proof of our main theorem, Theorem \ref{Main_theorem}. In Section \ref{Lemmas} we prove some lemmas that are used in earlier sections.

For simplicity we set $M = m = 0$ in the rest of the paper, but our proof can be modified to handle the massive case as well.

The following spaces of Bourgain-Klainerman-Machedon type are used throughout the paper. 
For $a,b \in \R$, we let $X^{a,b}_\pm$, $H^{a,b}$ and $\mathcal H^{a,b}$ be the completions of $\mathcal S(\R^{1+1})$ with respect to the norms
\begin{equation}\label{Spaces}
\begin{aligned}
  \norm{u}_{X^{a,b}_\pm} &= \bignorm{\angles{\xi}^a
\angles{\tau\pm \xi}^b \widetilde u(\tau,\xi)}_{L^2_{\tau,\xi}},
  \\
  \norm{u}_{H^{a,b}} &= \bignorm{\angles{\xi}^a \angles{\abs{\tau} -
\abs{\xi}}^b \widetilde u(\tau,\xi)}_{L^2_{\tau,\xi}},
  \\
  \norm{u}_{\mathcal H^{a,b}} &= \norm{u}_{H^{a,b}} + \norm{\partial_t u}_{H^{a-1,b}},
\end{aligned}
\end{equation}
where $\widetilde u(\tau,\xi)$ denotes the space-time Fourier transform of $u(t,x)$, and $\angles{\cdot} = 1 + \abs{\cdot}$. The restrictions of these spaces to a time slab
$$
  S_T = (0,T) \times \R
$$
are denoted $X^{a,b}_\pm(S_T)$ etc.; see \cite{Selberg:2006b} for more details about these spaces. We shall need the following basic estimates (see \cite{Selberg:2006b} for further references):

\begin{lemma}\label{X_lemma}
Let $1/2 < b \le 1$, $a \in \R$, $0 < T \le 1$ and $0 \le
\delta \le 1-b$. Then for all data $F \in X_\pm^{a,b-1+\delta}(S_T)$ and
$f \in H^a$, the Cauchy problem
$$
  (D_t \pm D_x) u = F(t,x) \quad \text{in $\,(0,T) \times \R$}, \qquad u(0,x) = f(x),
$$
has a unique solution $u \in X_\pm^{a,b}(S_T)$. Moreover,
$$
  \norm{u}_{X_\pm^{a,b}(S_T)} \le C \left( \norm{f}_{H^a}
  + T^\delta \norm{F}_{X_\pm^{a,b-1+\delta}(S_T)} \right),
$$
where $C$ only depends on $b$.
\end{lemma}
  
\begin{lemma}\label{H_lemma} \cite[Theorem 12]{Selberg:1999}.
Let $1/2 < b < 1$, $s \in \R$, $0 < T \le 1$ and $0 \le \delta
\le 1-b$. Then for all $F \in H^{a-1,b-1+\delta}(S_T)$, $f
\in H^a$ and $g \in H^{a-1}$, there exists a unique $u \in
\mathcal H^{a,b}(S_T)$ solving
$$
  \square u = F(t,x) \quad \text{in $\,(0,T) \times \R$}, \qquad u(0,x) = f(x),
  \qquad \partial_t u(0,x) = g(x).
$$
Moreover,
$$
  \norm{u}_{\mathcal H^{a,b}(S_T)} \le C \left( \norm{f}_{H^a} +
  \norm{g}_{H^{a-1}} + T^{\delta/2}
  \norm{F}_{H^{a-1,b-1+\delta}(S_T)} \right),
$$
where $C$ depends only on $b$.
\end{lemma}
We shall also need the fact that (see \cite{Selberg:2006b}) if $b > 1/2$, then
\begin{equation}\label{Basic_embedding}
  \norm{u(t)}_{H^a} \le C \norm{u}_{H^{a,b}(S_T)} \le C \norm{u}_{X_\pm^{a,b}(S_T)}
  \qquad \text{for $0 \le t \le T$},
\end{equation}
where $C$ only depends on $b$.

\section{Preliminaries}\label{Preliminaries}

We first observe that by propagation of higher regularity (see Remark \ref{Remark2} in Section \ref{Section1}), it suffices to prove Theorem \ref{Main_theorem} for $r < 1/2$. We therefore fix $s$ and $r$ satisfying
\begin{equation}\label{sr_cond}
  0 < s < \frac{1}{8},
  \qquad
  s + \sqrt{s^2+s} < r < \frac12,
\end{equation}
and set
\begin{equation}\label{b_def}
  b=\frac12+\varepsilon,
\end{equation}
where $\varepsilon > 0$ will be chosen sufficiently small, depending on $s$ and $r$; $b$ will be used as the second exponent in the spaces defined by \eqref{Spaces}.

Fix an arbitrary time $0 < T < \infty$. We shall split the interval $[0,T]$ into subintervals of length $\Delta T$, where the small time $\Delta T > 0$ remains to be chosen, and prove well-posedness on each subinterval successively. Let $M$ be the number of subintervals, so $M = T/\Delta T$.

Following Bourgain's idea, we split the data into low and high frequency parts:
$$
  \psi_0 = \psi_0^L + \psi_0^H,
$$
where $$(\psi_0^L)\FT(\xi) = \mathbb{1}_{\abs{\xi} \le N} \widehat \psi_0(\xi),$$ and the cut-off parameter $N \gg 1$ will be chosen later. Then
\begin{alignat}{2}
  \label{LowFreq}
  \norm{\psi_0^L}_{H^{-\zeta}} &\le C N^{s-\zeta}& \qquad &\text{for $\zeta \le s$},
  \\
  \label{HighFreq}
  \norm{\psi_0^H}_{H^{-\zeta}} &\le C N^{s-\zeta}& \qquad &\text{for $\zeta \ge s$},
\end{alignat}
where the constant $C$ depends on $\fixednorm{\psi_0}_{H^{-s}}$, but is independent of $N$. Note that we do not split the data for $\phi$, only for the spinor $\psi$.

Throughout the paper we use $C$ to denote a constant which can change from line to line; $C$ may depend on $s$, $r$, $T$ and $\varepsilon$, but it is always independent of $\Delta T$ and $N$.

We shall use $o(1)$ to denote a real-valued function $f(\varepsilon)$ such that $f(\varepsilon) \to 0$ as $\varepsilon \to 0$; these functions $f$ may depend implicitly on $s$ and $r$.

On each subinterval $[(n-1)\Delta T,n \Delta T]$ we split the solution of DKG into two parts:
\begin{equation}\label{induction_split}
  \psi = u_n + v_n, \qquad \phi = w_n + z_n \qquad \text{on} \qquad
  [(n-1)\Delta T,n \Delta T] \times \R,
\end{equation}
where $(u_n,w_n)$ solves DKG:
\begin{equation}\label{DKGuw}
\left\{
\begin{aligned}
  &\left(D_t + \alpha D_x\right)u_n  = w_n \beta u_n,
  \\
  &\square w_n = -\innerprod{\beta u_n}{u_n}_{\C^2},
\end{aligned}
\right.
\end{equation}
and hence the remaining part $(v_n,z_n) = (\psi-u_n,\phi-w_n)$ satisfies
\begin{equation}\label{DKGvz}
\left\{
\begin{aligned}
  &\left(D_t + \alpha D_x\right)v_n  = \phi\beta\psi - w_n \beta u_n
  \equiv z_n \beta v_n + z_n \beta u_n + w_n \beta v_n,
  \\
  &\square z_n = \innerprod{\beta u_n}{u_n}_{\C^2} - \innerprod{\beta \psi}{\psi}_{\C^2} \equiv -\innerprod{\beta v_n}{v_n}_{\C^2}-2\re \innerprod{\beta u_n}{v_n}_{\C^2}.
\end{aligned}
\right.
\end{equation}

For the first subinterval we take initial data
\begin{equation}\label{First_data}
\begin{alignedat}{3}
   u(0) &= \psi_0^L \in L^2,&
   \qquad w(0) &= \phi_0 \in H^r,&
   \qquad \partial_t w(0) &= \phi_1 \in H^{r-1},
   \\
   v(0) &= \psi_0^H \in H^{-s},&
   z(0) &= 0,&
   \partial_t z(0) &= 0.
\end{alignedat}
\end{equation}
Le $v^{(0)}$ denote the free evolution of $\psi_0^H$ by the Dirac equation:
\begin{equation}\label{v_free}
  (D_t + \alpha D_x ) v^{(0)} = 0, \qquad  v^{(0)}(0) = \psi_0^H \in H^{-s}.
\end{equation}
By Lemma \ref{X_lemma} and the estimate \eqref{HighFreq},
\begin{equation}\label{v:free:a}
  \bignorm{v_\pm^{(0)}}_{X_\pm^{-\zeta,b}(S_T)} \le C\norm{\psi_0^H}_{H^{-\zeta}} \le C N^{s-\zeta} \qquad \text{for $\zeta \ge s$}.
\end{equation}
Note that this holds on the whole time interval $[0,T]$, hence for any subinterval also.

For $(v_n,z_n)$ we then specify data, at time $t = (n-1) \Delta T$,
\begin{equation}\label{vz_data}
  v_n\bigl((n-1)\Delta T\bigr) = v^{(0)}\bigl((n-1)\Delta T\bigr),
  \qquad z_n\bigl((n-1)\Delta T\bigr) = \partial_t z_n\bigl((n-1)\Delta T\bigr) = 0.
\end{equation}
Writing $V_n$ for the inhomogeneous part of $v_n$ we then have
\begin{equation}\label{v_split}
  v_n = v^{(0)} + V_n \qquad \text{on} \qquad
  [(n-1)\Delta T,n \Delta T] \times \R.
\end{equation}
Of course, the free part $v^{(0)}$ just retains the regularity of the data. The inhomogeneous part $V_n(t)$, on the other hand, turns out to be smoother than the data, and in fact has $L^2$ regularity. This nonlinear smoothing effect is due to the null structure of the system, and it allows us to implement Bourgain's idea for the present problem, feeding $V_n$ back into the data for $u_{n+1}$, at time $n\Delta T$. We shall also feed $z_n$ into $w_{n+1}$.

Therefore, we specify the initial data for $(u_n,w_n)$ by the induction scheme
\begin{equation}\label{uw_data}
\begin{gathered}
  u_{n+1}(n\Delta T) = u_{n}(n\Delta T) + V_{n}(n\Delta T),
  \\
  w_{n+1}(n\Delta T) = w_{n}(n\Delta T) + z_{n}(n\Delta T),
  \quad \partial_t w_{n+1}(n\Delta T) =
  \partial_t w_{n}(n\Delta T) + \partial_t z_{n}(n\Delta T)
\end{gathered}
\end{equation}
for $1 \le n < M$. Recall that $M$ denotes the number of subintervals.

The main induction hypotheses will be:
\begin{gather}
  \label{u_n:data_bound}
  \norm{u_n\bigl((n-1)\Delta T\bigr)}_{L^2} \le A_n N^{s},
  \\
  \label{w_n:data_bound}
  \norm{w_n\bigl((n-1)\Delta T\bigr)}_{H^{r}} + \norm{\partial_t w_n\bigl((n-1)\Delta T\bigr)}_{H^{r-1}} \le B_n N^{2s},
\end{gather}
for all large $N$, where $A_n$ and $B_n$ are independent of $N$.
At the first induction step $n=1$, \eqref{u_n:data_bound} is satisfied in view of \eqref{LowFreq}, whereas the left side of \eqref{w_n:data_bound} does not depend on $N$ at all, so \eqref{w_n:data_bound} is trivially satisfied by choosing $B_1$ large enough. Note also that $A_1$ and $B_1$ do not depend $\varepsilon$, but this will not be the case for subsequent induction steps.

Thus, the underlying idea is simple, but to make it work in practice is not entirely trivial; we need quite strong estimates on $u$, $V$, $w$ and $z$ to close the induction. These estimates are derived in the next two sections. 

\section{Estimates for $u$ and $w$}\label{uw_estimates}

Here we describe the general induction step from time zero to time $\Delta T$, hence we drop the subscript $n$. Thus, we shall derive estimates for $(u,w)$ solving DKG:
\begin{equation}\label{DKGuw:2}
  \left(D_t + \alpha D_x\right)u  = w \beta u,
  \qquad
  \square w = -\innerprod{\beta u}{u}_{\C^2},
\end{equation}
with
\begin{equation}\label{uw:data_bounds}
  \norm{u(0)}_{L^2} \le AN^{s},
  \qquad
  \norm{w(0)}_{H^{r}} + \norm{\partial_t w(0)}_{H^{r-1}} \le BN^{2s},
\end{equation}
for all large $N$, where $A$ and $B$ are independent of $N$. In fact, the solution exists globally, by Theorem \ref{GWP_theorem}.

By the conservation of charge \eqref{conservation_of_charge},
\begin{equation}\label{u:cons}
  \norm{u(t)}_{L^2} = \norm{u(0)}_{L^2}
\end{equation}
for all $t$.

In order to derive sufficiently strong spacetime estimates for the solution, we need to exploit the null structure inherent in DKG; to reveal this structure we decompose the spinor $\psi$ using the eigenspace projections of the Dirac operator $\alpha D_x$, following \cite{Selberg:2006b, Selberg:2006e}. The symbol $\alpha \xi$ of $\alpha D_x$ has eigenvalues $\pm \xi$ and corresponding eigenspace projections
$$
  P_{\pm} = \frac12
  \begin{pmatrix}
    1 & \pm1  \\
    \pm1 & 1
  \end{pmatrix}.
$$
We then write
$$
  u = u_+ + u_-,
  \qquad \text{where} \qquad
  u_+ = P_+ u,
  \quad
  u_- = P_- u.
$$
Applying $P_{\pm}$ to the first equation in
(\ref{DKG1}), and using the identities $\alpha = P_+ -  P_-$, $P_{\pm}^2 = P_{\pm}$ and $P_+P_-=P_-P_+=0$, and the fact that $\beta$ is hermitian, \eqref{DKGuw:2} is rewritten as
\begin{equation}\label{DKGuw:3}
\left\{
\begin{aligned}
  &(D_t +  D_x ) u_+ = P_+(w \beta u) \equiv w\beta u_-,
  \\
  &(D_t -  D_x ) u_- = P_-(w \beta u) \equiv w\beta u_+,
  \\
  &\square w = -2\re \innerprod{\beta u_+}{u_-}_{\C^2},
\end{aligned}
\right.
\end{equation}
where we also used $\beta P_+ = P_- \beta$ and $\beta P_- = P_+ \beta$. The intrinsic null structure of DKG manifests itself through the different signs in the right hand side of the last equation; the same structure is in fact encoded in the first two equations, as becomes apparent via a duality argument; see \cite{Selberg:2006b, Selberg:2006e}.

We write
$$
  \norm{u} \equiv \norm{u_+}_{X_+^{0,b}(S_{\Delta T})} + \norm{u_-}_{X_-^{0,b}(S_{\Delta T})}.
$$

To estimate $w$ we shall use the following bilinear spacetime estimate for free waves in 1d, proved in \cite{Selberg:2006e}. A related estimate was proved by Bournaveas \cite[Lemma 1]{Bournaveas:2000}. The null structure, i.e., the difference of signs, is crucial here. In fact, the following is a 1d version of the null form estimate of Klainerman and Machedon \cite{Klainerman:1993}.

\begin{lemma}\label{WaveLemma} \cite{Selberg:2006e}.
Suppose $u,v$ solve
$$
\begin{alignedat}{2}
  &(D_t + D_x ) u = 0,& \qquad &u(0,x) = f(x),
  \\
  &(D_t - D_x ) v = 0,& \qquad &v(0,x) = g(x),
\end{alignedat}
$$
where $f,g \in L^2(\R)$. Then
$$
  \norm{uv}_{L^2(\R^{1+1})} \le \sqrt{2} \norm{f}_{L^2} \norm{g}_{L^2}.
$$
\end{lemma}

By the transfer principle (see \cite[Lemma 4]{Selberg:2006b}), Lemma \ref{WaveLemma} immediately implies:

\begin{corollary}\label{WaveCorollary} For any $b > 1/2$, the estimate
$$
  \norm{uv}_{L^2(\R^{1+1})} \le C \norm{u}_{X_+^{0,b}} \norm{v}_{X_-^{0,b}}
$$
holds, where $C$ depends only on $b$.
\end{corollary}

By Lemma \ref{H_lemma}, the induction hypothesis \eqref{uw:data_bounds} and the above corollary,
\begin{equation}\label{w:spacetime:a}
  \norm{w}_{H^{r,b}(S_{\Delta T})}
  \le
  C \left( \norm{w(0)}_{H^{r}} + \norm{\partial_t w(0)}_{H^{r-1}}
  +
  \norm{\square w}_{L^2(S_{\Delta T})}  \right)
  \le CBN^{2s} + C \norm{u}^2,
\end{equation}
since $\square w = -2\re\innerprod{\beta u_+}{u_-}$ by \eqref{DKGuw:3}.

To estimate $u$ we need the following lemma, proved in Section \ref{Lemmas}. Here there is again a crucial null structure due to the difference of signs between the left and right sides in the estimate.

\begin{lemma}\label{Lemma5}
Suppose $u$ solves
$$
  (D_t + D_x) u = \Phi \beta U, \qquad u(0) = u_0,
$$
where $\Phi$ is real-valued and $U$ is a 2-spinor. Assume $0 < r < 1/2$. Then with $b=1/2+\varepsilon$ and $\varepsilon > 0$ sufficiently small, depending on $r$, we have the estimate
$$
  \norm{u}_{X_+^{0,b}(S_{\Delta T})} \le C\norm{u_0}_{L^2} +  C (\Delta T)^{2r-7\varepsilon} \norm{\Phi}_{H^{r,b}(S_{\Delta T})} \norm{U}_{X_-^{0,b}(S_{\Delta T})}
$$
for all $0 < \Delta T \le T$, where $C$ only depends on $r$, $\varepsilon$ and $T$.
\end{lemma}

Applying this lemma to \eqref{DKGuw:3} and using \eqref{uw:data_bounds} and \eqref{w:spacetime:a}, we get
\begin{align*}
  \norm{u_+}_{X_+^{0,b}(S_{\Delta T})}
  &\le C \norm{u(0)}_{L^2}
  + C(\Delta T)^{2r-7\varepsilon} \norm{w}_{H^{r,b}(S_{\Delta T})}
  \norm{u_-}_{X_-^{0,b}(S_{\Delta T})}
  \\
  &\le CA N^s
  +  C (\Delta T)^{2r-7\varepsilon} B N^{2s}\norm{u}
  + C (\Delta T)^{2r-7\varepsilon} \norm{u}^3.
\end{align*}
Of course, the same estimate holds for $u_-$ in $X_-^{0,b}(S_{\Delta T})$, hence
\begin{equation}\label{u_est}
  \norm{u}
  \le CA N^s
  +  C (\Delta T)^{2r-7\varepsilon} B N^{2s}\norm{u}
  +  C (\Delta T)^{2r-7\varepsilon} \norm{u}^3,
\end{equation}
so if (with the same $C$)
\begin{equation}\label{contraction}
  C (\Delta T)^{2r-7\varepsilon} N^{2s} \left( 2B + 8C^2 A^2 \right)
  \le 1,
\end{equation}
then it follows by a boot-strap argument outlined below that (again with the same $C$)
\begin{equation}\label{u:spacetime}
  \norm{u}
  \le 2C A N^s.
\end{equation}
Then by \eqref{w:spacetime:a} we get also (now modifying $C$)
\begin{equation}\label{w:spacetime}
  \norm{w}_{H^{r,b}(S_{\Delta T})}
  \le C\left( B + A^2 \right) N^{2s}. 
\end{equation}

Motivated by the condition \eqref{contraction}, we now make the following choice:
\begin{equation}\label{time}
  \Delta T \sim N^{-(2s+\varepsilon)/(2r-7\varepsilon)}.
\end{equation}
Then $(\Delta T)^{2r-7\varepsilon} N^{2s} \sim N^{-\varepsilon}$
tends to zero as $N$ tends to infinity, hence \eqref{contraction} holds for $N$ large enough, provided $A$ and $B$ remain bounded throughout the induction; for the time being, however, \eqref{contraction} should be considered an induction hypothesis.

\begin{remark}\label{Remark3}
For the boot-strap argument referred to above we use the Picard iterates. Alternatively, one could use a continuity argument, but we want to avoid this since the continuity of, say, $\norm{u_+}_{X_+^{0,b}(S_{\Delta T})}$ as a function of $\Delta T$ is not obvious. The iterates $u^{(n)}$ are defined by setting $u^{(-1)} \equiv 0$ and then
$$
  (D_t \pm D_x ) u_\pm^{(n+1)} = w^{(n)} \beta u_\mp^{(n)}, \qquad u_\pm^{(n+1)}(0) = u_\pm(0),
  \qquad \text{for}\quad n = -1,0,1,\dots, 
$$
where $\square w^{(n)} = -2\re\biginnerprod{\beta u_+^{(n)}}{u_-^{(n)}}$ with the same data as for $w$.
Then \eqref{u_est} takes the form
$$
  y_{n+1} \le CA N^s
  +  C (\Delta T)^{2r-7\varepsilon} B N^{2s}y_n
  +  C (\Delta T)^{2r-7\varepsilon} y_n^3,
$$
for $n \ge 0$, where $y_n = \norm{u^{(n)}}$. By Lemma \ref{X_lemma} and \eqref{uw:data_bounds}, $y_0 \le CAN^s$, so if \eqref{contraction} holds, it now follows by induction that $y_n \le 2CAN^s$ for all $n \ge 1$. Since $\norm{u^{(n)}-u}$ tends to zero as $n$ tends to infinity (see \cite{Selberg:2006e}), this proves \eqref{u:spacetime}.
\end{remark}

\section{Estimates for $v,z$}\label{vz_estimates}

Again, we drop the subscript $n$, since we consider the general induction step.
We first reformulate \eqref{DKGvz} using the splitting $v=v_+ + v_-$, where $v_+ = P_+v$ and $v_- = P_-v$. Thus, we prove estimates on the time interval $[0,\Delta T]$ for $(v_+,v_-,z)$ solving
\begin{equation}\label{DKGvz:2}
\left\{
\begin{aligned}
  &(D_t +  D_x ) v_+  = z \beta v_- + z \beta u_- + w \beta v_-,
  \\
  &(D_t -  D_x ) v_-  = z \beta v_+ + z \beta u_+ + w \beta v_+,
  \\
  &\square z = -2\re\innerprod{\beta v_+}{v_-}_{\C^2}-2\re \innerprod{\beta u_+}{v_-}_{\C^2}
  -2\re \innerprod{\beta u_-}{v_+}_{\C^2},
\end{aligned}
\right.
\end{equation}
with $z(0) = \partial_t z(0) = 0$.
Here $u$ and $w$ are exactly as in the previous section. From now on we drop the subscript $\C^2$ on the inner product.

We write $v = v^{(0)} + V$,
where $v^{(0)}$ is the free part, which we assume satisfies (recall \eqref{v:free:a}), 
\begin{equation}\label{v:free}
  \bignorm{v^{(0)}}_{X_\pm^{-\zeta,b}(S_{\Delta T})}
  \le C N^{s-\zeta} \qquad \text{for $\zeta \ge s$}.
\end{equation}

Consider now the inhomogeneous part $V$ of $v$. We write
$$
  \norm{V_\pm} \equiv \norm{V_\pm}_{X_\pm^{0,b}(S_{\Delta T})},
  \qquad
  \norm{V} \equiv
  \norm{V_+} + \norm{V_-},
$$
and
$$
  V_\pm = V'_\pm + V''_\pm + V'''_\pm,
$$
where (see \eqref{DKGvz:2})
\begin{align*}
  (D_t \pm D_x) V'_\pm &= z\beta v_\mp,
  \\
  (D_t \pm D_x) V''_\pm &= w\beta v_\mp,
  \\
  (D_t \pm D_x) V'''_\pm &= z\beta u_\mp,
\end{align*}
with zero data at time zero.

We need the following null form estimate, proved in Section \ref{Lemmas}.

\begin{lemma}\label{Lemma2}
Let $b=1/2+\varepsilon$, where $0 < \varepsilon \le 1/2$, and assume $0 \le r \le b$. Then we have the 2-spinor estimate (note the different signs on the right)
$$
  \norm{\innerprod{\beta U}{V}}_{H^{-r,-b}} \le C \norm{U}_{X_+^{-r+2\varepsilon,b}} \norm{V}_{X_-^{0,1-b}}.
$$
Thus, by duality,
$$
  \norm{\Phi \beta U}_{X_-^{0,b-1}} \le C \norm{\Phi}_{H^{r,b}} \norm{U}_{X_+^{-r+2\varepsilon,b}}
$$
where $\Phi$ is real-valued and $U$ is a 2-spinor. Here $C$ depends only on $r$ and $\varepsilon$.
\end{lemma}

Using Lemmas \ref{X_lemma}, \ref{Lemma5} and \ref{Lemma2}, and the estimate \eqref{v:free}, we get
\begin{equation}\label{V:1}
\begin{aligned}
  \bignorm{V'_+}
  &\le C \norm{z \beta v^{(0)}_-}_{X_+^{0,b-1}(S_{\Delta T})}
  + C(\Delta T)^{2r-7\varepsilon}
  \norm{z}_{H^{r,b}(S_{\Delta T})} \norm{V}
  \\
  &\le C \norm{z}_{H^{r,b}(S_{\Delta T})}
  \norm{v_-^{(0)}}_{X_-^{-r+2\varepsilon,b}(S_{\Delta T})}
  + C(\Delta T)^{2r-7\varepsilon}
  \norm{z}_{H^{r,b}(S_{\Delta T})} \norm{V}
  \\
  &\le C \norm{z}_{H^{r,b}(S_{\Delta T})}
  \left[ N^{s-r+2\varepsilon}
  + (\Delta T)^{2r-7\varepsilon} \norm{V} \right].
\end{aligned}
\end{equation}
Note that we are justified in applying \eqref{v:free}, since $r > s$ by the assumption \eqref{sr_cond}.

Similarly, using \eqref{w:spacetime},
\begin{equation}\label{V:2}
\begin{aligned}
  \bignorm{V''_+}
  &\le C \norm{w}_{H^{r,b}(S_{\Delta T})}
  \left[ N^{s-r+2\varepsilon}
  + (\Delta T)^{2r-7\varepsilon} \norm{V} \right]
  \\
  &\le C \left( B + A^2 \right) N^{2s}
  \left[ N^{s-r+2\varepsilon}
  + (\Delta T)^{2r-7\varepsilon} \norm{V} \right].
\end{aligned}
\end{equation}

By Lemma \ref{Lemma5} and \eqref{u:spacetime},
\begin{equation}\label{V:3}
  \bignorm{V'''_+}
  \le C(\Delta T)^{2r-7\varepsilon} \norm{z}_{H^{r,b}(S_{\Delta T})}
  \norm{u}
  \le CA N^s(\Delta T)^{2r-7\varepsilon} \norm{z}_{H^{r,b}(S_{\Delta T})}.
\end{equation}

Next, we estimate $z$. Using \eqref{DKGvz:2} and the fact that $z$ has zero data at time zero, we write
\begin{equation}\label{z:spacetime:a}
  \norm{z}_{H^{r,b}(S_{\Delta T})}
  \le 2\left( I_1 + I_2 + I_3 + I_4 + I_5 \right) + (\cdots),
\end{equation}
where
\begin{align*}
  I_1 &= \norm{\square^{-1} \biginnerprod{\beta v^{(0)}_+}{v^{(0)}_-}}_{H^{r,b}(S_{\Delta T})},
  \\
  I_2 &= \norm{\square^{-1} \biginnerprod{\beta v^{(0)}_+}{V_-}}_{H^{r,b}(S_{\Delta T})},
  \\
  I_3 &= \norm{\square^{-1} \biginnerprod{\beta v^{(0)}_+}{u_-}}_{H^{r,b}(S_{\Delta T})},
  \\
  I_4 &= \norm{\square^{-1} \biginnerprod{\beta V_+}{V_-}}_{H^{r,b}(S_{\Delta T})},
  \\
  I_5 &= \norm{\square^{-1} \biginnerprod{\beta V_+}{u_-}}_{H^{r,b}(S_{\Delta T})},
\end{align*}
and $(\cdots)$ indicates similar terms for which we have the same estimates as for the $I_j$'s. Here we use the notation $\square^{-1} F$ for the solution of $\square Z = F$ with zero data at time zero.

To estimate $I_1$ we need the following, proved in Section \ref{Lemmas}.

\begin{lemma}\label{Lemma6}
Assume $r \le 1/2$ and $b=1/2+\varepsilon$, where $0 < \varepsilon \le 1/2$. Then
$$
  \norm{\innerprod{\beta U}{V}}_{H^{r-1,b-1}} \le C \norm{U}_{X_+^{-1/4+\varepsilon,b}} \norm{V}_{X_-^{-1/4+\varepsilon,b}},
$$
where $C$ depends only on $r$ and $\varepsilon$.
\end{lemma}

By Lemmas \ref{X_lemma}, \ref{Lemma6} and the estimate \eqref{v:free},
\begin{equation}\label{I:1}
\begin{aligned}
  I_1 &\le C\norm{\innerprod{\beta v^{(0)}_+}{v^{(0)}_-}}_{H^{r-1,b-1}(S_{\Delta T})}
  \\
  &\le C\bignorm{v_+^{(0)}}_{X_+^{-1/4+\varepsilon,b}(S_{\Delta T})}
  \bignorm{v_-^{(0)}}_{X_-^{-1/4+\varepsilon,b}(S_{\Delta T})}
  \\
  &\le C N^{2(s-1/4+\varepsilon)}.
\end{aligned}
\end{equation}

For $I_2$ and $I_3$ we need the following variation of Lemma \ref{Lemma6}.

\begin{lemma}\label{Lemma7}
Assume $r < 1/2$ and $b=1/2+\varepsilon$, where $0 < \varepsilon \le 1/2$. Then
$$
  \norm{\innerprod{\beta U}{V}}_{H^{r-1,b-1}} \le C \norm{U}_{X_+^{-1/2+\varepsilon,b}} \norm{V}_{X_-^{0,b}},
$$
where $C$ depends only on $r$ and $\varepsilon$.
\end{lemma}

Combining this with Lemma \ref{X_lemma}, \eqref{v:free} and \eqref{u:spacetime}, we get
\begin{align}
  \label{I:2}
  I_2 &\le C\bignorm{v_+^{(0)}}_{X_+^{-1/2+\varepsilon,b}(S_{\Delta T})} \norm{V}
  \le CN^{s-1/2+\varepsilon} \norm{V},
  \\
  \label{I:3}
  I_3 &\le C\bignorm{v_+^{(0)}}_{X_+^{-1/2+\varepsilon,b}(S_{\Delta T})} \norm{u} \le C AN^{2s-1/2+\varepsilon}.
\end{align}

Finally, to estimate $I_4$ and $I_5$ we need the following. Note that in this estimate the signs do not matter.

\begin{lemma}\label{Lemma8}
Suppose $u$ solves
$$
  \square u = \innerprod{\beta U}{V}, \qquad u \init = \partial_t u \init = 0,
$$
where $U,V$ are 2-spinors. Assume $r < 1/2$. Then with $b=1/2+\varepsilon$ and $\varepsilon > 0$ sufficiently small, depending on $r$, we have
$$
  \norm{u}_{H^{r,b}(S_{\Delta T})} \le C (\Delta T)^{3/4-2\varepsilon} \norm{U}_{X_+^{0,b}(S_{\Delta T})} \norm{V}_{X_\pm^{0,b}(S_{\Delta T})}
$$
for all $0 < \Delta T \le T$, where $C$ only depends on $r$, $\varepsilon$ and $T$.
\end{lemma}

This lemma gives
\begin{equation}\label{I:4}
  I_4 \le C (\Delta T)^{3/4-2\varepsilon} \norm{V}^2,
\end{equation}
and, using \eqref{u:spacetime},
\begin{equation}\label{I:5}
  I_5 \le C (\Delta T)^{3/4-2\varepsilon} \norm{V} \norm{u}
  \le CA N^s (\Delta T)^{3/4-2\varepsilon} \norm{V}.
\end{equation}

From \eqref{z:spacetime:a}--\eqref{I:5} we conclude that
$$
  \norm{z}_{H^{r,b}(S_{\Delta T})}
  \le C AN^{2s-1/2+2\varepsilon}
  + CAN^s \left[(\Delta T)^{3/4-2\varepsilon} + N^{-1/2+\varepsilon}\right] \norm{V} 
  +  C (\Delta T)^{3/4-2\varepsilon} \norm{V}^2.
$$
However, $N^{-1/2+\varepsilon}$ is negligible compared to $(\Delta T)^{3/4-2\varepsilon}$ if $\varepsilon > 0$ is sufficiently small and $N$ sufficiently large, since by \eqref{time}, $\Delta T \sim N^{-s/r+o(1)}$, and $r > 2s$ by \eqref{sr_cond}. Hence,
\begin{equation}\label{z:spacetime}
  \norm{z}_{H^{r,b}(S_{\Delta T})}
  \le C AN^{2s-1/2+2\varepsilon}
  + CAN^s (\Delta T)^{3/4-2\varepsilon} \norm{V} 
  +  C (\Delta T)^{3/4-2\varepsilon} \norm{V}^2.
\end{equation}
Combining this with \eqref{V:1}, \eqref{V:2} and \eqref{V:3}, we conclude:
\begin{equation}\label{V:boot_strap:1}
\begin{aligned}
  \norm{V} 
  &\le
  C \left(B + A^2\right) N^{2s} \left[ N^{s-r+2\varepsilon} + N^{s-1/2+2\varepsilon}(\Delta T)^{2r-7\varepsilon} \right]
  \\
  &\quad
  + C\left(B + A^2\right) N^{2s} \left[ (\Delta T)^{2r-7\varepsilon} + (\Delta T)^{3/4+2r-9\varepsilon}
  + N^{-r+2\varepsilon}(\Delta T)^{3/4-2\varepsilon} \right]
  \norm{V}
  \\
  &\quad
  + CAN^{s}  (\Delta T)^{3/4-2\varepsilon}\left[ (\Delta T)^{2r-7\varepsilon}
  + N^{-r+2\varepsilon}
  \right]
  \norm{V}^2
  \\
  &\quad
  + C (\Delta T)^{3/4+2r-9\varepsilon} \norm{V}^3.
\end{aligned}
\end{equation}
We want to simplify the right hand side. Consider first the sum inside the parentheses in the first line. Clearly, the term $N^{s-r+2\varepsilon}$ dominates, since $r < 1/2$. Now look at terms in parentheses in the third line. The term $(\Delta T)^{2r-7\varepsilon}$ dominates, since $(\Delta T)^{2r-7\varepsilon} \sim N^{-2s+2\varepsilon}$, by \eqref{time}, and since $r > 2s$; the same reasoning applies for the second line.

Thus,
\begin{equation}\label{V:boot_strap:2}
\begin{aligned}
  \norm{V} 
  &\le
  C \left(B + A^2\right) N^{3s-r+2\varepsilon}
  + C\left(B + A^2\right) N^{2s} (\Delta T)^{2r-7\varepsilon}
  \norm{V}
  \\
  &\quad
  + CAN^{s}  (\Delta T)^{3/4+2r-9\varepsilon}  \norm{V}^2
  + C (\Delta T)^{3/4+2r-9\varepsilon} \norm{V}^3.
\end{aligned}
\end{equation}
But in view of the induction hypothesis \eqref{contraction}, the term which is linear in $\norm{V}$ can be moved to the left hand side, yielding
\begin{equation}\label{V:boot_strap:3}
  \norm{V} 
  \le
  C \left(B + A^2\right) N^{3s-r+2\varepsilon}
  + CAN^{s}  (\Delta T)^{3/4+2r-9\varepsilon}  \norm{V}^2
  + C (\Delta T)^{3/4+2r-9\varepsilon} \norm{V}^3.
\end{equation}
A boot-strap argument described below then gives (with the same $C$ as in \eqref{V:boot_strap:3})
\begin{equation}\label{V:final}
  \norm{V} \le 2\Gamma,
  \quad \text{where}
  \quad
  \Gamma = C \left(B + A^2\right) N^{3s-r+2\varepsilon},
\end{equation}
provided that (still with the same $C$) the following induction hypothesis holds:
\begin{equation}\label{contraction2}
  C (\Delta T)^{3/4+2r-9\varepsilon} \left[ 4AN^{s}\Gamma + 8 \Gamma^2 \right] \le 1.
\end{equation}

\begin{remark}\label{Remark4}
The above estimates imply local well-posedness of \eqref{DKGvz:2} for data $v(0) \in H^{-s}$ satisfying $\norm{v(0)}_{H^{-\zeta}} \le CN^{s-\zeta}$ for $\zeta \ge s$, and zero data for $z$, with existence up to the time $\Delta T > 0$ determined by the conditions \eqref{contraction} and \eqref{contraction2}, by a standard argument using the iterates. To define the iterates we write $v_\pm^{(n)} = v_\pm^{(0)} + V_\pm^{(n)}$ for $n = 1,2,\dots$, where $v^{(0)}$ is the free part and $V_\pm^{(n)}$ for $n = 1,2,\dots$ are determined by the scheme
$$
  (D_t \pm D_x) V^{(n+1)}_\pm = z^{(n)} \beta v_\mp^{(n)} + z^{(n)}\beta u_\mp  + w \beta v_\mp^{(n)}, \qquad  V^{(n+1)}_\pm(0) = 0,
  \qquad  \text{for $n \ge 0$},
$$
where
$$
  \square z^{(n)} = - 2\re\innerprod{\beta v_+^{(n)}}{v_-^{(n)}}
  - 2\re\innerprod{\beta u_+}{v_-^{(n)}}
  - 2\re\innerprod{\beta v_+^{(n)}}{u_-}
$$
with zero initial data. Then our estimates imply that $V_\pm^{(n)}$ converges to $V_\pm$ in $X_\pm^{0,b}(S_{\Delta T})$, as $n$ tends to infinity. The boot-strap procedure referred to above is included in this iteration argument; in this setting, \eqref{V:boot_strap:3} becomes
\begin{equation}
  y_{n+1} 
  \le
  \Gamma
  + CAN^{s}  (\Delta T)^{3/4+2r-9\varepsilon} y_n^2
  + C (\Delta T)^{3/4+2r-9\varepsilon} y_n^3,
  \qquad \text{for $n \ge 0$},
\end{equation}
where $y_0 = 0$ and $y_n = \norm{V^{(n)}}$ for $n \ge 1$. So if \eqref{contraction2} holds, then by induction we get $y_n \le 2\Gamma$ for all $n$.
\end{remark}

From \eqref{z:spacetime} and \eqref{V:final} we conclude that
$$
  \norm{z}_{H^{r,b}(S_{\Delta T})}
  \le
  C AN^{2s-1/2+2\varepsilon}
  + C (\Delta T)^{3/4-2\varepsilon}  \left( AN^s \Gamma 
  + \Gamma^2 \right).
$$
But since $r > 2s$,
\begin{equation}\label{gamma_est}
  \Gamma^2 \le C (B+A^2) N^s \Gamma,
\end{equation} 
so we get
\begin{equation}\label{z:final}
  \norm{z}_{H^{r,b}(S_{\Delta T})}
  \le
  C AN^{2s-1/2+2\varepsilon}
  + C (B+A^2)^2 N^{4s-r-3s/4r + o(1)},
\end{equation}
where we used that $\Delta T \sim N^{-s/r+o(1)}$, by \eqref{time}.

\section{Conclusion of the proof}\label{Conclusion}

We now apply the estimates proved in the last two sections to the induction scheme described in Section \ref{Preliminaries}. Recall that on the $n$-th subinterval $[(n-1)\Delta T,\Delta T]$, where $\Delta T$ satisfies \eqref{time}, $(u_n,w_n)$ solves \eqref{DKGuw} and $(v_n,z_n)$ solves \eqref{DKGvz}; the data prescribed at time $(n-1)\Delta T$ are given by \eqref{vz_data} for $(v_n,z_n)$ and by the induction scheme \eqref{uw_data} for $(u_n,w_n)$. The main induction hypotheses are \eqref{u_n:data_bound} and \eqref{w_n:data_bound}, involving the constants $A_n$ and $B_n$, which must be independent of the large parameter $N$. In addition, the boot-strap conditions \eqref{contraction} and \eqref{contraction2} must be satisfied at each step of the induction. Recall also that $M$ denotes the number of induction steps, so
\begin{equation}\label{M}
  M = \frac{T}{\Delta T} \sim N^{s/r+o(1)},
\end{equation}
where we used \eqref{time}.

Note that the boot-strap conditions \eqref{contraction} and \eqref{contraction2} reduce to, using \eqref{time} and \eqref{gamma_est},
\begin{equation}\label{contraction3}
  CN^{-\varepsilon}(B_n+A_n^2) \le 1
  \qquad \text{and} \qquad
  CN^{2s-r-(3s)/(4r) + o(1)} ( B_n + A_n^2 )^2 \le 1.
\end{equation}

The crucial point now is to prove that $A_n$ and $B_n$ do not grow indefinitely, but remain bounded for $1 \le n \le M$. Otherwise it would be impossible to satisfy \eqref{contraction3} with $N$ chosen independently of $n$.

By \eqref{uw_data}, the conservation of charge \eqref{conservation_of_charge}, the induction hypothesis \eqref{u_n:data_bound}, the embedding \eqref{Basic_embedding} and the estimate \eqref{V:final}, assuming \eqref{contraction3} holds,
\begin{equation}\label{S5:A}
\begin{aligned}
  \norm{u_{n+1}(n\Delta T)}_{L^2}
  &\le \norm{u_{n}(n\Delta T)}_{L^2} + \norm{V_{n}(n\Delta T)}_{L^2}
  \\
  &= \norm{u_{n}\bigl((n-1)\Delta T\bigr)}_{L^2} + \norm{V_{n}(n\Delta T)}_{L^2}
  \\
  &\le A_n N^s + \norm{V_{n}(n\Delta T)}_{L^2}
  \\
  &\le A_n N^s + C \left(B_n + A_n^2\right) N^{3s-r+2\varepsilon}.
\end{aligned}
\end{equation}
Therefore,
\begin{equation}\label{A_induction}
   A_{n+1} \le A_n +  C \left(B_n + A_n^2\right) N^{2s-r+2\varepsilon}.
\end{equation}

It will be convenient to use the notation
$$
  \norm{\phi[t]}_{H^r} \equiv \norm{\phi(t)}_{H^r} + \norm{\partial_t \phi(t)}_{H^{r-1}}.
$$
By \eqref{uw_data},
\begin{equation}\label{w_naive}
  \norm{w_{n+1}[n\Delta T]}_{H^{r}}
  \le
  \norm{w_n[n\Delta T]}_{H^{r}}
  + \norm{z_n[n\Delta T]}_{H^{r}},
\end{equation}
but the problem is that we have no conservation law for $w$. Implementing Bourgain's idea for DKG is therefore not so straightforward, and requires an additional idea.

We first split $w_n$ into its free part $w_n^{(0)}$ and its inhomogeneous part $W_n$, so
\begin{equation}\label{w_splitting}
  w_n = w_n^{(0)} + W_n,
\end{equation}
where $\square W_n = \square w_n$ and $W_n$ has zero data at time $(n-1)\Delta T$. The inhomogeneous part is quite favorable. In fact, applying the estimate \eqref{phi_apriori} to $W_n$ on the interval $[(n-1)\Delta T,\Delta T]$, and using \eqref{u_n:data_bound}, we get, assuming \eqref{contraction3} holds,
\begin{equation}\label{W_add}
  \norm{W_n[n\Delta T]}_{H^{r}}
  \le C (\Delta T) \norm{u_n\bigl((n-1)\Delta T\bigr)}_{L^2}^2
  \le C (\Delta T) A_n^2 N^{2s}.
\end{equation}
By \eqref{Basic_embedding} and \eqref{z:final}, assuming \eqref{contraction3} holds,
\begin{equation}\label{z_add}
\begin{aligned}
  \norm{z_n[n\Delta T]}_{H^{r}}
  &\le C \norm{z_n}_{H^{r,b}([(n-1)\Delta T,n\Delta T] \times \R)}
  \\
  &\le C A_n N^{2s-1/2+2\varepsilon}
  + C (B_n+A_n^2)^2 N^{4s-r-3s/4r + o(1)}.
\end{aligned}
\end{equation}
That leaves us with the free part $w_n^{(0)}$. Certainly, by the energy inequality and \eqref{w_n:data_bound},
$$
  \norm{w_n^{(0)}[n\Delta T]}_{H^{r}} \le C \norm{w_n^{(0)}[(n-1)\Delta T]}_{H^{r}} \le C B_n N^{2s},
$$
but this naive estimate is useless, since it gives at best $B_n \sim C^n$. However, things are not as bad as they may seem at a first glance. The important point is to keep accurately track of what is added to the free part at the end of each induction step. In fact, by induction, referring to the scheme \eqref{uw_data}, we have
\begin{equation}\label{cascade}
  w_n^{(0)} = w_1^{(0)} + \tilde w_2^{(0)} + \tilde w_3^{(0)} + \dots + \tilde w_n^{(0)},
\end{equation}
for $n \ge 1$, where
\begin{align*}
  \square \tilde w_n^{(0)} &= 0,
  \\
  \tilde w_n^{(0)}\bigl( (n-1)\Delta T \bigr) &=
  W_{n-1}\bigl( (n-1)\Delta T \bigr)
  + z_{n-1}\bigl( (n-1)\Delta T \bigr),
  \\
  \partial_t \tilde w_n^{(0)}\bigl( (n-1)\Delta T \bigr) &=
  \partial_t W_{n-1}\bigl( (n-1)\Delta T \bigr)
  + \partial_t z_{n-1}\bigl( (n-1)\Delta T \bigr).
\end{align*}
Thus, we express $w_n^{(0)}$ as a cascade of free waves. By the energy inequality for the wave equation, and the estimates \eqref{W_add} and \eqref{z_add},
\begin{equation}\label{w_tilde}
  \bignorm{\tilde w_{n+1}^{(0)}[t]}_{H^r} \le
  C (\Delta T) A_n^2 N^{2s}
  + C A_n N^{2s-1/2+2\varepsilon}
  + C (B_n+A_n^2)^2 N^{4s-r-3s/4r + o(1)},
\end{equation}
in the entire time interval $0 \le t \le T$.

We now replace the induction hypothesis \eqref{w_n:data_bound} by the stronger condition
\begin{equation}\label{w:free_bound}
  \sup_{0 \le t \le T} \norm{w_n^{(0)}[t]}_{H^r} \le B_n N^{2s}.
\end{equation}
In fact, by the energy inequality this condition is equivalent to \eqref{w_n:data_bound}, up to multiplication by a constant depending only on $T$.

Since $w_{n+1}^{(0)} = w_n^{(0)} + \tilde w_{n+1}^{(0)}$, we have
$$
  \bignorm{w_{n+1}^{(0)}[t]}_{H^r} \le \norm{w_n^{(0)}[t]}_{H^r}
  + \bignorm{\tilde w_{n+1}^{(0)}[t]}_{H^r},
$$
for all $0 \le t \le T$, and we conclude from \eqref{w_tilde} that
\begin{equation}\label{B_induction}
  B_{n+1}
  \le B_n
  + C (\Delta T) A_n^2
  + C A_n N^{-1/2+2\varepsilon}
  + C (B_n+A_n^2)^2 N^{2s-r-3s/4r + o(1)}.
\end{equation}

We claim that if $\varepsilon > 0$ is chosen small enough, and then $N$ large enough, depending on $\varepsilon$, then for $1 \le n \le M$,
\begin{equation}\label{AB_bound}
  A_n \le R \equiv 2A_1
  \qquad \text{and} \qquad
  B_n \le S \equiv 2B_1 + 4CT A_1^2
\end{equation}
with the same $C$ that appears in the second term of the right hand side of \eqref{B_induction}.

We proceed by induction. Assume \eqref{AB_bound} holds for $1 \le n < m$, for some $m \le M$.
Then \eqref{contraction3} reduces to
\begin{equation}\label{C1}
  CS N^{-\varepsilon} \le 1,
\end{equation}
and
\begin{equation}\label{C2}
  C S^2 N^{2s-r-(3r)/(4s) + o(1)} \le 1,
\end{equation}
for $n < m$. Since $r > 2s$, we can ensure that the exponent of $N$ in \eqref{C2} is negative by choosing $\varepsilon > 0$ small enough. Then $C$, and hence $S$, may grow, but this is not a problem, since we can now choose $N$ as large as we need, depending on $\varepsilon$, to ensure that \eqref{C1} and \eqref{C2} are satisfied. Then \eqref{A_induction} and \eqref{B_induction} are also satisfied for $n < m$, so by the assumption that \eqref{AB_bound} holds for $n < m$, we get
\begin{align*}
  A_{n+1} &\le A_1 + n CS N^{2s-r+o(1)},
  \\
  B_{n+1} &\le B_1 + n \left[ CRN^{-1/2+o(1)}+ CS^2 N^{2s-r-3s/(4r)+o(1)}
  + CR^2 \Delta T \right],
\end{align*}
for $n < m$. Thus, \eqref{AB_bound} will hold also for $A_m$ and $B_m$ provided that
\begin{gather*}
  (m-1) CS N^{2s-r+o(1)}
  \le A_1,
  \\
  (m-1) \left[ CR  N^{-1/2+o(1)} + CS^2 N^{2s-r-3s/(4r)+o(1)} + CR^2 (\Delta T) \right]
  \le B_1 + 4CTA_1^2.
\end{gather*}
But $m \le M \le C N^{s/r+o(1)}$, by \eqref{M}, so it suffices to have
\begin{gather}
  \label{Cond1}
  CS N^{2s-r+s/r+o(1)} \le A_1,
  \\
  \label{Cond2}
  CR N^{-1/2+s/r+o(1)} \le B_1/2,
  \\
  \label{Cond3}
  CS^2  N^{2s-r-3s/(4r)+s/r+o(1)} \le B_1/2,
  \\
  \label{Cond5}
  C T R^2 \le 4CTA_1^2,
\end{gather}
where to get \eqref{Cond5} we used the fact that
$m \Delta T \le M \Delta T = T$, by \eqref{M}.
Note that \eqref{Cond5} holds with equality since $R=2A_1$. To satisfy \eqref{Cond1}--\eqref{Cond3} it suffices to have
\begin{equation}\label{Conditions}
  2s-r+\frac{s}{r} < 0,
  \qquad
  -\frac12 +\frac{s}{r} < 0,
  \qquad
  2s-r-\frac{3s}{4r}+\frac{s}{r} < 0.
\end{equation}
Indeed, if these inequalities hold, then we first choose $\varepsilon > 0$ so small that the exponents of $N$ in \eqref{Cond1}--\eqref{Cond3} are also negative. Then the constants $C$ (hence also $S$) may grow, since they depend on $\varepsilon$. But by subsequently choosing $N$ large enough we can ensure that \eqref{Cond1}--\eqref{Cond3} hold.

The first inequality in \eqref{Conditions} is equivalent to $r^2-2sr-s > 0$, i.e.,
$r > s + \sqrt{s^2+s}$,
which holds by assumption \eqref{sr_cond}. The last two inequalities in \eqref{Conditions} are weaker than the first one.
 
Thus, \eqref{AB_bound} holds for $n=1,\dots,M$, and the proof of Theorem \ref{Main_theorem} is complete.

\section{Proof of Lemmas}\label{Lemmas}

\subsection{Proof of Lemma \ref{Lemma2}}

This estimate is a variation on an estimate proved in \cite{Selberg:2006e}, and we use the same method of proof as in that paper. First, by Plancherel's theorem, the estimate is equivalent to
\begin{equation}\label{I_def}
  I = \norm{ \int_{\R^{1+1}} \frac{ F(\lambda,\eta)
  G(\lambda-\tau,\eta-\xi)d\lambda \, d\eta}
  {\angles{\xi}^{r}\angles{\eta}^{-r+2\varepsilon} \angles{A}^{b}
  \angles{B_+}^{b}
  \angles{C_-}^{1-b}}}_{L^2_{\tau,\xi}}
  \le C \norm{F}_{L^2} \norm{G}_{L^2},
\end{equation}
for arbitrary $F, G \in L^2(\R^{1+1})$, where
$$
 A=\abs{\tau}-\abs{\xi}, \qquad
 B_+=\lambda+\eta, \qquad
 C_-=(\lambda-\tau)-(\eta-\xi).
$$

The null structure, i.e., the fact that there are two different signs, is crucial, since it allows us to use the following estimate, proved in \cite[Lemma 1]{Selberg:2006e}:
\begin{equation}\label{algebraic}
 \min \left( \fixedabs{\eta}, \fixedabs{\eta-\xi} \right)
 \le \frac32 \max\left( \abs{A}, \abs{B_+}, \abs{C_-} \right).
\end{equation}
In fact, this estimate is the 1d version of \cite[Lemma 7]{Selberg:2006b}.

Thus, we can reduce to proving \eqref{I_def} with $I$ replaced by one of the following:
\begin{align*}
  I_1 &= \norm{ \int_{\R^{1+1}} \frac{ F(\lambda,\eta)
  G(\lambda-\tau,\eta-\xi)d\lambda \, d\eta}
  {\angles{\xi}^{r}\angles{\eta}^{1-b-r+2\varepsilon} \angles{A}^{\theta_1}
  \angles{B_+}^{\theta_2}
  \angles{C_-}^{\theta_3}}}_{L^2_{\tau,\xi}},
  \\
   I_2 &= \norm{ \int_{\R^{1+1}} \frac{ F(\lambda,\eta)
  G(\lambda-\tau,\eta-\xi)d\lambda \, d\eta}
  {\angles{\xi}^{r}\angles{\eta}^{-r+2\varepsilon} \angles{\eta-\xi}^{1-b} \angles{A}^{\theta_1}
  \angles{B_+}^{\theta_2}
  \angles{C_-}^{\theta_3}}}_{L^2_{\tau,\xi}},
\end{align*}
where the $\theta_j$ are nonnegative and $\theta_1+\theta_2+\theta_3 > 1/2$.
Then we can apply the following product law for the spaces $H^{s,b}$.

\begin{theorem}\label{H_products} \cite[Proposition 10]{Selberg:1999}.
Suppose $s_1,s_2,s_3 \in \R$, \,$\theta_1, \theta_2, \theta_3 \ge 0$, and
$$
  \theta_1+\theta_2+\theta_3 > \frac12,
  \qquad
  s_1+s_2+s_3 > \frac12,
  \qquad s_1+s_2 \ge 0,
  \qquad s_1+s_3 \ge 0,
  \qquad s_2+s_3 \ge 0.
$$
Then (in one space dimension)
$$
  \norm{uv}_{H^{-s_1,-\theta_1}} \le C \norm{u}_{H^{s_2,\theta_2}} \norm{v}_{H^{s_3,\theta_3}},
$$
where $C$ depends on the $s_j$ and the $\,\theta_j$.
\end{theorem}

\begin{remark} In \cite{Selberg:1999} only the case $s_1,s_2,s_3 \ge 0$ is proved, but the general case reduces to this by Leibniz' rule (the triangle inequality in Fourier space).
\end{remark}

Applying this to $I_1$ and $I_2$, and recalling that $b=1/2+\varepsilon$, we see that it suffices to have $0 \le r \le 1/2+\varepsilon = b$ and $0 < \varepsilon \le 1/2$. This completes the proof of Lemma \ref{Lemma2}.

\subsection{Proof of Lemma \ref{Lemma7}} Proceeding as in the proof of Lemma \ref{Lemma2}, we reduce to proving \eqref{I_def}, but with $I$ replaced by one of the following:
\begin{align*}
  I_1 &= \norm{ \int_{\R^{1+1}} \frac{ F(\lambda,\eta)
  G(\lambda-\tau,\eta-\xi)d\lambda \, d\eta}
  {\angles{\xi}^{1-r}\angles{\eta}^{1-b-1/2+\varepsilon} \angles{A}^{\theta_1}
  \angles{B_+}^{\theta_2}
  \angles{C_-}^{\theta_3}}}_{L^2_{\tau,\xi}},
  \\
   I_2 &= \norm{ \int_{\R^{1+1}} \frac{ F(\lambda,\eta)
  G(\lambda-\tau,\eta-\xi)d\lambda \, d\eta}
  {\angles{\xi}^{1-r}\angles{\eta}^{-1/2+\varepsilon} \angles{\eta-\xi}^{1-b} \angles{A}^{\theta_1}
  \angles{B_+}^{\theta_2}
  \angles{C_-}^{\theta_3}}}_{L^2_{\tau,\xi}},
\end{align*}
where the $\theta_j$ are nonnegative and $\theta_1+\theta_2+\theta_3 > 1/2$. Applying Theorem \ref{H_products} and recalling that $b=1/2+\varepsilon$ with $0 < \varepsilon \le 1/2$, we get the sufficent condition $r < 1/2$.

\subsection{Proof of Lemma \ref{Lemma6}} This reduces to \eqref{I_def} with $I$ replaced by
$$ 
  I = \norm{ \int_{\R^{1+1}} \frac{ F(\lambda,\eta)
  G(\lambda-\tau,\eta-\xi)d\lambda \, d\eta}
  {\angles{\xi}^{1-r}\angles{\eta}^{1-b-1/4+\varepsilon}\angles{\eta-\xi}^{-1/4+\varepsilon} \angles{A}^{\theta_1}
  \angles{B_+}^{\theta_2}
  \angles{C_-}^{\theta_3}}}_{L^2_{\tau,\xi}},
$$
where the $\theta_j$ are nonnegative and $\theta_1+\theta_2+\theta_3 > 1/2$. Since  $b=1/2+\varepsilon$ with $0 < \varepsilon \le 1/2$, Theorem \ref{H_products} yields the sufficent condition $r \le 1/2$.

\subsection{Proof of Lemma \ref{Lemma8}}

Applying Lemma \ref{H_lemma}, H\"older's inequality in time, the Sobolev product law \eqref{Sobolev_product}, and the embedding \eqref{Basic_embedding}, we get for all $r < 1/2$,
\begin{align*}
  \norm{u}_{H^{r,b}(S_{\Delta T})}
  &\le C (\Delta T)^{(1-b)/2} 
  \norm{\innerprod{\beta U}{V}}_{{L_t^2H^{r-1}(S_{\Delta T})}}
  \\
  &\le C (\Delta T)^{(1-b)/2} (\Delta T)^{1/2} 
  \norm{\innerprod{\beta U}{V}}_{{L_t^\infty H^{r-1}(S_{\Delta T})}}
  \\
  &\le C (\Delta T)^{(1-b)/2} (\Delta T)^{1/2} 
  \norm{U}_{{L_t^\infty L_x^2(S_{\Delta T})}}
  \norm{V}_{{L_t^\infty L_x^2(S_{\Delta T})}}
  \\
  &\le C (\Delta T)^{3/4-\varepsilon} \norm{U}_{X_+^{0,b}(S_{\Delta T})} \norm{V}_{X_\pm^{0,b}(S_{\Delta T})}.
\end{align*}

\subsection{Proof of Lemma \ref{Lemma5}}

It suffices to prove this for $u_0 = 0$, in view of Lemma \ref{X_lemma}. Then by Lemma \ref{X_lemma}, H\"older's inequality in time, the product law \eqref{Sobolev_product}, and the embedding \eqref{Basic_embedding},
\begin{align*}
  \norm{u}_{X_+^{0,b}(S_{\Delta T})}
  &\le C (\Delta T)^{1/2-\varepsilon} (\Delta T)^{1/2}
  \norm{\Phi\beta U}_{L_t^\infty L_x^2(S_{\Delta T})}
  \\
  &\le C  (\Delta T)^{1-\varepsilon} \norm{\Phi}_{H^{1/2+\varepsilon,b}(S_{\Delta T})} \norm{U}_{X_-^{0,b}(S_{\Delta T})}.
\end{align*}
On the other hand, by Lemmas \ref{X_lemma} and \ref{Lemma2},
$$
  \norm{u}_{X_+^{0,b}(S_{\Delta T})}
  \le C \norm{\Phi\beta U}_{X_+^{0,b-1}(S_{\Delta T})}
  \le C \norm{\Phi}_{H^{2\varepsilon,b}(S_{\Delta T})} \norm{U}_{X_-^{0,b}(S_{\Delta T})}.
$$
Interpolation gives, for $0 \le \theta \le 1$,
$$
  \norm{u}_{X_+^{0,b}(S_{\Delta T})}
  \le C  (\Delta T)^{\theta(1-\varepsilon)} \norm{\Phi}_{H^{(1-\theta)2\varepsilon+\theta(1/2+\varepsilon),b}(S_{\Delta T})} \norm{U}_{X_-^{0,b}(S_{\Delta T})}.
$$
If $0 < r < 1/2$, then taking $\theta = 2r-6\varepsilon$ with $\varepsilon$ small enough, we get the desired result, since then
$(1-\theta)2\varepsilon+\theta(1/2+\varepsilon) \le \theta/2 + 3\varepsilon = r$,
and $\theta(1-\varepsilon) \ge 2r-7\varepsilon$.

\bibliographystyle{amsplain}
\bibliography{/Users/sselberg/Mathematics/Bibliography/mybibliography}

\end{document}